\documentclass{amsart}
\usepackage{amsmath,amscd,amssymb,latexsym,amsfonts,epsfig}

\renewcommand{\proof}{\par\noindent{\it Proof.\ \ }}
\def\qed{\ifmmode\square\else\nolinebreak\hfill
$\Box$\fi\par\vskip12pt}

\input xy
\xyoption{all}
\usepackage[all,cmtip]{xy}
\newtheorem{thm}{Theorem}[section]
\newtheorem{prop}[thm]{Proposition}
\newtheorem{lem}[thm]{Lemma}
\newtheorem{cor}[thm]{Corollary}

\theoremstyle{definition}
\newtheorem{rem}[thm]{Remark}

\def\Aut{{\rm Aut}}

\def\Ga{\Gamma}

\newcommand{\Z}{\mathbb Z}

\newcommand{\C}{\mathbb C}
\newcommand{\ra}{\rightarrow}

\newcommand{\<}{\langle}
\renewcommand{\>}{\rangle}

\begin{document}
\pagestyle{plain}

\begin{titlepage}

\begin{center}
\author{ Yu Qing Chen, Henry H. Glover and Craig A. Jensen}
\end{center}

\address{Department of Mathematics and Statistics, Wright State
University, Dayton, OH 45435-0001} \email{yuqing.chen@wright.edu}

\address{Department of Mathematics, Ohio State University, Columbus,
OH 43210-1174} \email{glover@math.ohio-state.edu}

\address{Department of Mathematics, University of New Orleans, New Orleans,
LA 70148} \email{jensen@math.uno.edu}

\begin{abstract}
In this paper we first show that many braid groups of low genus
surfaces have their centers as direct factors. We then give a
description of centralizers and normalizers of prime order elements
in pure mapping class groups of surfaces with spherical quotients
using automorphism groups of fundamental groups of the quotient
surfaces. As an application, we use these to show that the
$p$-primary part of the Farrell cohomology groups of certain mapping
class groups are elementary abelian groups. At the end we compute
the $p$-primary part of the Farrell cohomology of a few pure mapping
class groups.
\end{abstract}

\title{ The center of some braid groups and the Farrell cohomology of
certain pure mapping class groups}


\maketitle

\end{titlepage}

\section{Introduction}

Braid groups were first introduced by Artin \cite{A1,A2,A3} for the
purpose of studying knots and links. A geometric braid on $n$
strings can be defined as follow. Let $E={\mathbb E}^2\times I,$
where ${\mathbb E}^2$ is the Euclidean plane and $I$ is the unit
interval $[0,1]$. Let $P_i=(i,0,1)\in{\mathbb E}^2\times\{1\}$ and
$P^{'}_i=(i,0,0)\in{\mathbb E}^2\times\{0\}$ for $i=1,2,\cdots,n.$
An $n$-string $s=\{s_1,s_2,\cdots,s_n\}$ in $E$ is a set of maps
$s_i:[0,1]\rightarrow E,$ $i=1,2,\cdots,n,$ satisfying
\begin{itemize}
\item[(I)] $s_i(0)=P_i$ and $s_i(1)=P^{'}_{\sigma(i)},$
$i=1,2,\cdots,n,$ where $\sigma$ is a permutation of the set
$\{1,2,\cdots,n\}$, and
\item[(II)] each plane ${\mathbb E}^2\times\{t\}$ in $E$, $0\leqslant
t\leqslant 1,$ intersects the strings $\bigcup_{i=1}^ns_i([0,1])$ in
exactly $n$ points.
\end{itemize}
Two $n$-strings $s$ and $s^{'}$ given by the same permutation
$\sigma$ of $\{1,2,\cdots,n\}$ are said to be equivalent if there is
a homotopy $H=(H_1,H_2,\cdots,H_n)$, where $H_i,$ $i=1,2,\cdots,n$,
are maps
$$H_i:[0,1]\times[0,1]\rightarrow E,$$ such that
$$\left\{ \begin{array}{l}
              H_i(u,0)=s_i(u), \\
              H_i(u,1)=s^{'}_i(u),\\
              H_i(0,v)=P_i,\\
              H_i(1,v)=P^{'}_{\sigma(i)}
         \end{array}
\right.$$ and such that for every $v\in [0,1],$
$H_v=\{H_1(\cdot,v),H_2(\cdot,v) ,\cdots,H_n(\cdot,v)\}$ is an
$n$-string. Then a geometric braid on $n$ strings is an equivalence
class of $n$-strings. Geometric braids can be presented pictorially
by two horizontal lines with the points $P_i$'s on the top line and
$P^{'}_i$'s on the bottom line and strings going from $P_i$'s to
$P^{'}_i$'s. The set of braids on $n$ strings is denoted by $B_n.$
We can define a product of two braids in $B_n$ using these pictures.
Let $b_1$ and $b_2$ be two geometric braids on $n$ strings. The
picture of $b=b_1b_2$ is obtained as follow: placing the picture of
$b_1$ above $b_2$ and identifying the bottom line of $b_1$ with the
top line of $b_2$ in such a way that the points
$P^{'}_1,P^{'}_2\cdots,P^{'}_n$ in $b_1$ match $P_1,P_2,\cdots,P_n$
in $b_2.$ Then removing the common line between $b_1$ and $b_2$. It
is easy to check that this product is well defined. Therefore $B_n$
is a group which is called Artin braid group on $n$ strings. A braid
with the trivial permutation $\sigma$ is called a pure braid. The
set of pure braids forms a normal subgroup of $B_n$ and is denoted
by $P_n$. It is called the pure Artin braid group. There is a exact
sequence \begin{eqnarray}\label{ex1}1\ra P_n\ra B_n\ra\Sigma_n\ra
1,\end{eqnarray} where $\Sigma_n$ is the symmetry group on $n$
symbols. Given a subgroup $G$ of $\Sigma_n$, we define $G$-Artin
braid group $B^G_n$ on $n$ strings to be the subgroup of $B_n$ such
that for all braids in $B^G_n$, the permutations $\sigma$ in (I) are
from $G$. Clearly, $P_n\leqslant B^G_n\leqslant B_n$ and
$B^{\Sigma_n}_n=B_n$ and $B^{\{1\}^n}_n=P_n$, where $\{1\}^k$ is the
trivial subgroup of $\Sigma_k$. The group $B^G_n$ is simply the
pre-image of $G$ in $B_n$ of the exact sequence (\ref{ex1}) and one
has an exact sequence
$$1\ra P_n\ra B^G_n\ra G\ra 1.$$

Another family of groups that we are interested in are the mapping
class groups. They are closely related to braid groups. Let $S_g^n$
denote the orientable surface of genus $g$ with $n$ punctures. The
mapping class group $\Gamma_g^n$ of $S_g^n$ is defined to be the
isotopy classes of orientation preserving homeomorphisms from
$S_g^n$ to itself. It is the discrete group of connected components
of the group of orientation preserving homeomorphisms of $S_g^n$.
For other equivalent definitions of $\Gamma_g^n$ and the properties
of $\Gamma_g^n$ with $n=0,$ we refer the reader to the survey paper
of Mislin \cite{M}. Each $\Gamma_g^n$ contains a normal subgroup
consisting of isotopy classes of orientation preserving
homeomorphisms of $S_g^n$ which point-wise fix the punctures. This
subgroup is called pure mapping class group of $S_g^n$ and is
denoted by $P\Gamma_g^n$. Similar to the braid groups, we have the exact
sequence \begin{eqnarray}\label{ex2} 1\ra
P\Ga_g^n\ra\Ga_g^n\ra\Sigma_n\ra 1.\end{eqnarray} For each subgroup
$G$ of $\Sigma_n$, we can also define the $G$-mapping class group of
$S^n_g$, $\Ga_g^{n,G}$, to be the pre-image of $G$ in $\Ga_g^n$ of
the exact sequence (\ref{ex2}). Then one has
$P\Ga_g^n\leqslant\Ga_g^{n,G}\leqslant\Ga_g^n$, and
$\Ga_g^{n,\Sigma_n}=\Ga_g^n$ and $\Ga_g^{n,\{1\}^n}=P\Ga_g^n$, and
the exact sequence $$1\ra P\Ga_g^n\ra\Ga_g^{n,G}\ra G\ra 1.$$

In \cite{Lu1}, Lu proved that $P\Ga_g^n$ has periodic cohomology
with period 2 for all $g\geqslant 1$ and $n\geqslant 1$, i.e.
$\widehat{H}^{i+2}(P\Ga_g^n,\Z)\cong\widehat{H}^{i}(P\Ga_g^n,\Z)$.
She also computed the $p$-primary part of  the Farrell cohomology of
punctured pure mapping class group
$\widehat{H}^*(P\Ga_{n(p-1)/2}^m,\Z)_{(p)}$ for some small values of
$n$ ($n\leqslant 3$) as well as $\widehat{H}^*(P\Ga_{g}^n,\Z)$ for
small values of $g$ ($g\leqslant 3$)\cite{Lu1,Lu2,Lu3}.

In this paper, we first show in Theorems~\ref{split}, \ref{splitS}
and \ref{splitT} that for many braid groups of $S^1_0$(2-sphere with
1 puncture), $S^0_0$(2-sphere with no punctures) and $S^0_1$(torus
with no punctures) (see the definition of these groups in Section
2), their centers are actually direct factors of these groups. We
then give in Theorem~\ref{ppart} an algebraic description of
centralizers and normalizers of elements of order $p$ in $P\Ga_g^m$
when the quotient surfaces are punctured 2-spheres. Combining the
splitting of braid groups and the alternative description of
normalizers, we are able to obtain structural information about the
$p$-primary part of the Farrell cohomology of $P\Ga_{n(p-1)/2}^m$
for certain values of $n$ and $m$ in Theorem~\ref{ele}. At the end
of the of the paper in Theorem~\ref{explicit}, we give an explicit
computation of ${\widehat H}^i(P\Gamma_{n(p-1)/2}^{n+1},{\mathbb
Z})_{(p)}$, ${\widehat H}^i(P\Gamma_{n(p-1)/2}^{n+2},{\mathbb
Z})_{(p)}$, ${\widehat H}^i(P\Gamma_{n}^{2n+1},{\mathbb Z})_{(2)}$
and ${\widehat H}^i(P\Gamma_{n}^{2n+2},{\mathbb Z})_{(2)}$.

The rest of the paper is organized as follow. In section 2 we review
some well known facts about braid groups and mapping class groups.
In section 3  we prove that the braid groups of $S^1_0$, $S^0_0$ and
$S^0_1$ are the direct product of their centers with other braid
groups or mapping class groups. In section $4$, we compute
centralizers and normalizers of elements of prime order in pure
mapping class groups, provided that the quotient surfaces are
punctured 2-spheres. In the last section, we study the $p$-primary
part of the Farrell cohomology of pure mapping class groups and
compute the $p$-primary part of Farrell cohomology mentioned above.

\section{Braid groups and mapping class groups of surfaces}
This section contains the descriptions of braid groups and mapping
class groups of surfaces which are different from the ones given in
the introduction. The relations between pure braid groups and pure
mapping class groups of the punctured 2-spheres are discussed at the
end of this section.

The use of configuration spaces to describe braid groups seems to be
suggested by Fox. This approach enables one to define braid groups of any
topological space, and the braid groups defined in section 1 are
simply the braid groups of Euclidean plane.

Let $X$ be a topological space. Its $n$-fold
configuration space $F_n(X)$ for $n\geq 0$ is defined by
$$F_n(X)=\{(x_1,x_2,\cdots,x_n)\in\prod_{i=1}^{n}X~|~x_i\ne x_j\mbox{
for all } i\ne j\}$$ with the convention that when $n=0$, the space
$F_0(X)$ consists of a single point and when $n=1,$ the space
$F_1(X)=X$ for any space $X$. Each subgroup $G$ of $\Sigma_n$ has a
natural free action on  $F_n(X)$ by permuting the coordinates. We
denote by $F_n^G(X)$ the quotient space $F_n(X)/G$, i.e.
$$F_n^G(X)=\{(x_1,x_2,\cdots,x_n)^G~|~x_i\ne x_j \mbox{
for all } i\ne j\},$$ where $(x_1,x_2,\cdots,x_n)^G$ is the
$G$-orbit of $(x_1,x_2,\cdots,x_n)$ in $F_n(X)$. Then the $G$-braid
group of $X$ on $n$ strings is defined to be the fundamental group
$\pi_1(F_n^G(X))$ of $F_n^G(X)$. When $G=\{1\}^n$,
$F_n^{\{1\}^n}(X)=F_n(X)$ and $\pi_1(F_n(X))$ is called the pure braid
group of $X$ on $n$ strings. When $G=\Sigma_n$,
$\pi_1(F_n^{\Sigma_n}(X))$ is called full braid group of $X$ on $n$
strings. Obviously, the natural projection from $F_n(X)$ to
$F_n^G(X)=F_n(X)/G$ is a regular covering with deck transformation
group $G$. Thus we have an exact sequence of groups
$$1\rightarrow\pi_1(F_n(X))\rightarrow\pi_1(F_n^G(X))\rightarrow G
\rightarrow 1.$$ If the space $X$ is the Euclidean plane, or
equivalently, the 1-punctured 2-sphere $S_0^1$, then the braid group
$\pi_1(F_n^G(S_0^1))$ of $S_0^1$ is precisely the $G$-Artin braid
group $B^G_n$. The subgroup $\pi_1(F_n(S_0^1))$ is the pure Artin
braid group $P_n.$ For the connection between the two definitions of
the Artin braid groups we refer the reader to \cite{Bi}. For a fixed
set of distinguished points $Q_m=\{x_1,x_2,\cdots,x_m\}$ of a
topological space $X,$ we define
$$F_{m,n}(X)=F_n(X-Q_m)$$ and for any subgroup $G$ of $\Sigma_n$ we
define $$F^G_{m,n}(X)=F^G_n(X-Q_m).$$ Note that the topological type
of the configuration space $F_{m,n}(X)$ is independent of the choice
of the particular points in the set $Q_m$, since one may always find
an isotopy of $X$ which deforms one such set to another. Let
$n=s+t$, $G=G_s\times G_t$ for some subgroups $G_s$ in $\Sigma_s$
and $G_t$ in $\Sigma_t$.  The following theorem, which is a slight
generalization of a version in \cite{FN}, describes a link among
various configuration spaces $F^G_{m,n}(X)$ for a given space $X$.

{\thm\label{FN} \cite[Fadell and Neuwirth]{FN} Let $\pi:
F^{G_s\times F_t}_{m,s+t}(X)\rightarrow F^{G_s}_{m,s}(X)$ be the
projection defined by
$$\pi((x_1,x_2,\cdots,x_{s+t})^{G_s\times G_t})=(x_1,x_2,\cdots,x_s)^{G_s}.$$
 Then $\pi$ is a locally trivial fibration
with fibre $F^{G_t}_{m+s,t}(X).$} {\rem The fibration is usually not
 globally trivial.} \vspace{1.5mm}

We will use this theorem to establish a connection between
$\pi_1(F^G_n(S_0^1))$ and $\pi_1(F^G_n(S_0^0))$.

In \cite{A1}, Artin gave an algebraic description of the Artin braid
groups which is now called the Artin representation.
{\thm\label{AlgArtin}\cite[Artin]{A1} Let  $Aut(F_n)$ denote the group of automorphisms
 of the free group with generators $x_1,x_2,\cdots x_n.$  The Artin braid group
\begin{align*}
B_n\cong ~&\{\gamma\in {\rm Aut}(F_n)~|~\gamma(x_i)\mbox{ is
conjugate to }
x_{\sigma(i)} \mbox{ for some permutation } \sigma \mbox{ of }\\
&\{1,2,\cdots,n\}\mbox{ and for all } i=1,2,\cdots,n \mbox{ and }
\gamma(x_1x_2\cdots x_n)=x_1x_2\cdots x_n \}.
\end{align*}

}\vspace{1.5mm}

 It follows easily from Theorem~\ref{AlgArtin} that for any subgroup
 $G$ of $\Sigma_n$, the $G$-Artin braid group
\begin{align*}
B_n^G\cong~&\{\gamma\in {\rm Aut}(F_n)~|~\gamma(x_i)\mbox{ is
conjugate to }x_{\sigma(i)} \mbox{ for some }\sigma\in G\\&\mbox{
and for all }i=1,2,\cdots,n\mbox{ and }\gamma(x_1x_2\cdots
x_n)=x_1x_2\cdots x_n\}.
\end{align*}
The centers of $G$-Artin braid groups are infinite cyclic groups.
Geometrically, these groups are generated by the braids with a full
twist. If we use the Artin representation to describe these centers,
they are precisely the intersection of the braid groups with the
inner automorphism groups of those free groups.

Mapping class groups $\Ga_g^{n,G}$ introduced in the previous
section can also be given an algebraic description using
automorphisms of fundamental groups of the surfaces $S_g^n$ as in
\cite{DF,Z}. The fundamental group of $S_g^n$ admits a presentation
\begin{align*}
\pi_1(S_g^n)=\langle x_1,x_2,\cdots,x_n,a_1,b_1,a_2,b_2,
\cdots,a_g,b_g ~|~x_1x_2\cdots x_n\prod_{i=1}^n[a_i,b_i]=1\rangle,
\end{align*}
where $x_1,$ $x_2,$ $\cdots,$ $x_n$ are the loops representing those
punctures which are oriented in a consistent manner. For a subgroup
$G$ of $\Sigma_n$,  we define
\begin{align*}
{\rm Aut}^G_+(\pi_1(S_g^n))=~&\{\gamma\in {\rm
Aut}(\pi_1(S_g^n))~|~\gamma(x_i) \mbox{ is conjugate to
}x_{\sigma(i)}\\&\mbox{ for some }\sigma\in G \mbox{ and for all }
i=1,2,\cdots, n\}.
\end{align*}
Then the $G$-mapping class group of $S^n_g$
\begin{align*}
\Gamma_g^{n,G}\cong~&{\rm Aut}^G_+(\pi_1(S_g^n))/{\rm
Inn}(\pi_1(S_g^n)),
\end{align*} where ${\rm Inn}(\pi_1(S_g^n))$ is the inner
automorphism group of $\pi_1(S_g^n)$.

If we set $g=0$, by comparing  the algebraic descriptions of $B^G_n$
and $\Ga^{n,G}_0$, it is obvious that
$B^G_n/Z(B^G_n)\cong\Ga_0^{n+1,G\times\{1\}}$, where $Z(B^G_n)$ is
the center of $B^G_n$.

We end this section with the following proposition that summarizes
the known relation between the braid groups of $S_0^1$ and $S_0^0$
and the mapping class groups of the punctured 2-sphere.  Some
version of this can be found in \cite{Bi,FN}.

{\prop

\begin{itemize}
\item[(1)] The center $Z(B^G_n)$ of $B^G_n$ is an infinite cyclic group and
$B^G_n/Z(B^G_n)\cong \Gamma_0^{n+1,G\times\{1\}}$;
\item[(2)] the center $Z(\pi_1(F^G_n(S_0)))$ of $\pi_1(F^G_n(S_0))$
is a cyclic group of order $2$ and
$\pi_1(F^G_n(S_0))/Z(\pi_1(F^G_n(S_0)))\cong \Gamma_0^{n,G}$.
\end{itemize}}

\section{The center of braid groups}

In this section we will show that under some mild conditions, the
center of many braid groups of
$S_0^1$, $S_0^0$, $S_1^0$ are direct factors of these groups. Let us
recall a few more well known facts about the braid groups. The
configuration space $F_n(S_0^1)$, which we used to define the pure
Artin group $P_n$, is an Eilenberg-MacLane space $K(P_n,1)$ of $P_n$
since $\pi_i(F_n(S_0^1))=0$ for $i\geqslant 2$ (see \cite{Bi}). Also
since $F_n(S_0^1)$ is a finite cover of $F^G_n(S_0^1)$ for any
subgroup $G$ of $\Sigma_n$, one has $\pi_i(F^G_n(S_0^1))=0$ for
$i\geqslant 2$. This implies that $B^G_n$ is torsion free because it
has finite cohomological dimension. In fact, If $S$ is a surface
with $\pi_2(S)=1$, i.e. $S\ne S^0_0$, the space $F^G_n(S)$ is an
Eilenberg-MacLane space of $\pi_1(F_n^G(S))$ for any $n$ and any
subgroup $G$ of $\Sigma_n$. Hence these braid groups are of finite
cohomological dimension and are torsion free. Also from  the
geometric definition of $P_n$, there is an exact sequence
$$1\rightarrow F_{n-1}\rightarrow P_n\rightarrow P_{n-1}\rightarrow
1,$$ where the projection $P_n\rightarrow P_{n-1}$ is given by
removing a string from $P_n$ and the kernel of this projection is
isomorphic to the free group $F_{n-1}$ on $n-1$ generators. The
above exact sequence splits. Therefore $$P_n\cong F_{n-1}\rtimes
P_{n-1},$$ where the action of $P_{n-1}$ on $F_{n-1}$ is given by
the Artin representation. For example, $P_3\cong F_2\rtimes{\mathbb
Z}$ and the action of ${\mathbb Z}$ on $F_2$ is given by an inner
automorphism of $F_2.$

We now prove a simple lemma concerning the configuration spaces of
topological groups.

{\lem If $K$ is a topological group and $G$ is a subgroup of
$\Sigma_n$ such that $G=H\times\{1\}$, where $H$ is a subgroup of
$\Sigma_{n-1}$, then
$$F_n^G(K)\cong F_{n-1}^H(K-\{{\bf 1}\})\times K,$$
where ${\bf 1}$ is the identity element of $K$.}

\proof We define
\begin{eqnarray*}
&&u:F^G_n(K)\rightarrow F^H_{n-1}(K-\{{\bf 1}\})\times K,\\
&&u((k_1,k_2,\cdots,k_{n-1},k_n)^G)=
((k_1k_n^{-1},k_2k_n^{-1},\cdots,k_{n-1}k_n^{-1})^H,k_n)
\end{eqnarray*}
and
\begin{eqnarray*}
&&v:F^H_{n-1}(K-\{{\bf 1}\})\times K\rightarrow F^G_n(K),\\
&&v((k_1,k_2,\cdots,k_{n-1})^H,k_n)=(k_1k_n,k_2k_n,\cdots,k_{n-1}k_n,k_n)^G.
\end{eqnarray*}
It is easy to verify that $u$ and $v$ are continuous and inverse to
each other. \qed

{\rem Comparing this Lemma with Theorem 2.1, we can see that the
fibration in Theorem 2.1 with $m=0$ and $s=1$ is usually not a
trivial fibration, while Lemma 3.1 provides a trivial one when the
space has a group structure.}

{\thm\label{split} Let $G$ be a subgroup of $\Sigma_n$ such that
$G=H\times\{1\}^2$, where $H$ is a subgroup of $\Sigma_{n-2}$. The
$G$-Artin braid group
$$B^G_n\cong \pi_1(F^H_{n-2}(S^3_0))\times\Z
\cong\Gamma_0^{n+1,G\times\{1\}}\times{\mathbb Z}$$ for $n\geqslant
2,$ where $G\times\{1\}$ is a subgroup of
$\Sigma_{n+1}$.}\vspace{1.5mm}

\proof  We use the complex number plane ${\mathbb C}$ for the
Euclidean plane ${\mathbb E}^2$ and denote ${\mathbb C}-\{0\}$ by
${\mathbb C}^*$, the non-zero complex numbers. Then we have
\begin{align*}
B^G_n=~&\pi_1(F^G_n({\mathbb C}))\\
     \cong~&\pi_1(F^{H\times\{1\}}_{n-1}({\mathbb C}^*)\times{\mathbb C})\\
     =~&\pi_1(F^{H\times\{1\}}_{n-1}({\mathbb C}^*))\\
     \cong~&\pi_1(F^H_{n-2}({\mathbb C}^*-\{1\})\times{\mathbb C}^*)\\
     \cong~&\pi_1(F^H_{n-2}({\mathbb C}^*-\{1\}))\times\pi_1({\mathbb
  C}^*)\\
     =~&\pi_1(F^H_{n-2}(S^3_0))\times\Z
\end{align*}
Here we used Lemma 3.1 twice, once for the additive group of
${\mathbb C}$ and once for the multiplicative group ${\mathbb C}^*.$
The factor $\pi_1({\mathbb C}^*)\cong{\mathbb Z}$ is obviously
contained in the center $Z(B^G_n)$ of $B^G_n$. In order to show that
it is actually the full center, we only need that
$Z(B^G_n)\cong{\mathbb Z}$. By statement (1) of Proposition 2.4,
we have
$$B^G_n\cong \Gamma_0^{n+1,G\times\{1\}}\times{\mathbb Z}.$$ \qed
{\cor\label{puresplit}
$$P_n\cong P\Gamma_0^{n+1}\times{\mathbb Z}.$$}
{\rem We can also
use the fact that $P_n\cong F_{n-1}\rtimes P_{n-1}$ to prove
Corollary~\ref{puresplit} by induction.}
{\rem Note that in Theorem 3.3, we have $G=H\times\{1\}^2$, a subgroup
of $\Sigma_n$ that fixes at least two elements. If $G$ fixes at most one
element, i.e. $G$ does not fix any elements or $G=H\times\{1\}$ with $H$
a subgroup of $\Sigma_{n-1}$ which does not fix any elements, then the
braid group $B_n^{G}$ does not split.}
{\cor The space
$F^H_{n-3}(S^3_0)$ is an Eilenberg-MacLane space of
$\Gamma_0^{n,H\times\{1\}^3}$ for any subgroup $H$ of
$\Sigma_{n-3}$. Therefore $\Gamma_0^{n,H\times\{1\}^3}$ has finite
cohomological dimension and is torsion free.} {\rem The corollary
recovers the well known fact that $P\Gamma_0^4\cong F_2,$ the free
group on two generators. Furthermore, combining this fact with
Corollary 3.4, we obtain that $P_3\cong F_2\times{\mathbb Z}.$ This
description of $P_3$ is simpler than the one given before Lemma
3.1.} {\rem When $H$ is a subgroup of $\Sigma_{n-2}$, the group
$\Ga_0^{n,H\times\{1\}^2}$ may contain torsion, and hence is not of
finite cohomological dimension.}\vspace{1.5mm}

We now use Lemma 3.1 to give a similar decomposition of the
$G$-braid group of the 2-sphere. {\thm\label{splitS} Let $G$ be a
subgroup of $\Sigma_{n}$ such that $G=H\times\{1\}^3$ for some
subgroup $H$ of $\Sigma_{n-3}$. The $G$-braid group of the
$2$-sphere
$$\pi_1(F^{G}_n(S_0^0))\cong \Gamma_0^{n,G}
\times{\mathbb Z}/2{\mathbb Z}$$ for $n\geqslant 3.$}\vspace{1.5mm}

\proof If we apply Theorem~\ref{FN} with $m=0,$ $s=1$, $t=n-1$ and
$X=S_0^0,$ we obtain a fibration $$\pi: F^{G}_n(S_0^0)\rightarrow
S^0_0$$ with fibre
$F^{H\times\{1\}^2}_{1,n-1}(S_0^0)=F^{H\times\{1\}^2}_{n-1}(S_0^1).$
Thus we have an exact sequence of the homotopy groups of the
fibration
\begin{eqnarray*}1=\pi_2(F^{H\times\{1\}^2}_{n-1}(S_0^1))\rightarrow\pi_2
(F^G_n(S_0^0))\rightarrow
\pi_2(S_0^0)\rightarrow\\\pi_1(F^{H\times\{1\}^2}_{n-1}(S_0^1))
\rightarrow\pi_1(F^G_n(S_0^0))
\rightarrow\pi_1(S_0^0)=1&.\end{eqnarray*} Since
$\pi_2(S_0^0)\cong{\mathbb Z}$ and, by Theorem~\ref{split},
$\pi_1(F^{H\times\{1\}^2}_{n-1}(S_0^1))\cong
\Gamma_0^{n,G}\times{\mathbb Z},$ the above exact sequence is
equivalent to
$$1\rightarrow\pi_2(F^G_n(S_0^0))\stackrel{i}{\rightarrow}{\mathbb
Z}\stackrel{j} {\rightarrow}\Gamma_0^{n,G}\times{\mathbb
Z}\stackrel{k}{\rightarrow}\pi_1(F^G_n(S_0^0) )\rightarrow 1.$$ The
surjection $k$ maps the center of $\Gamma_0^{n,G}\times{\mathbb Z}$
onto the center of $\pi_1(F_n(S_0^0))$ which is isomorphic to
${\mathbb Z}/2{\mathbb Z}.$ Therefore we must have that
$\pi_2(F^G_n(S_0^0))=1$ and $j({\mathbb Z})$ is a subgroup of index
2 of the center of $\Gamma_0^{n,G}\times{\mathbb Z}.$ From this one
can easily see that
$$\pi_1(F^G_n(S_0^0))\cong\Gamma_0^{n,G}\times{\mathbb Z}/2 {\mathbb Z}.$$
\qed

The last application of Lemma 3.1 is to prove the splitting of some braid
groups of torus.

{\thm\label{splitT} Let $G$ be a subgroup of $\Sigma_{n}$ such that
$G=H\times\{1\}$ for some subgroup $H$ of $\Sigma_{n-1}$. The
$G$-braid group of the torus
$$\pi_1(F^{G}_n(S_1^0))\cong \pi_1(F^{H}_{n-1}(S_1^1))
\times{\mathbb Z}^2$$ for $n\geqslant 1.$}\vspace{1.5mm}

\proof Let $S=\{x\in\C^*|~|x|=1\}$ be the unit circle in the complex number
plane. Torus $S_1^0$ is homeomorphic to $S\times S$, which is a topological
group. By Lemma 3.1,
\begin{align*}
\pi_1(F^{G}_n(S_1^0))=&\pi_1(F^{G}_n(S\times S))\\
                     \cong &\pi_1(F^{H}_{n-1}((S\times S)-\{(1,1)\})\times{(S\times S))}\\
                     \cong &\pi_1(F^{H}_{n-1}(S_1^1))\times\pi_1{(S\times S)}\\
                     \cong &\pi_1(F^{H}_{n-1}(S_1^1))\times{\mathbb Z}^2
\end{align*} for $n\geqslant 1.$

\section{The centralizer and normalizer of elements of prime order in mapping
class groups}

In this section we discuss the centralizer and normalizer of elements of
prime order in mapping class groups $\Ga_g^n$. Detailed discussion
of conjugacy classes of prime order finite subgroup of $P\Ga_g^n$
can be found in \cite{Lu1,Lu2}. Let us first recall two theorems
from \cite{BH}.

{\thm\label{BH1}\cite[Theorem 1]{BH} Let $(p,S_g^m, S_h^n)$ be a
regular covering space, either branched or unbranched, with a finite
group of covering transformations and a finite number of branch
points. Let the covering transformations leave each branch point
fixed. In the case of a branched covering, assume that $S_g^m$ is
not homeomorphic to the closed sphere or closed torus. Let
$g:S_g^m\ra S_g^m$ be a fiber-preserving homeomorphism of $S_g^m$
which is isotopic to the identity map. Then $g$ is fiber-isotopic to
the identity map.}\vspace{1.5mm}

The term ``fiber-isotopic" means the intermediate homeomorphisms
between $g$ and the identity map in the isotopy are all
fiber-preserving.

{\thm\label{BH3}\cite[Theorem 3]{BH} Let $t, h\in\Ga_g^0$ such that
$t^k=1$ and $hth^{-1}\in\<t\>$. There are homeomorphisms $\tilde{t}$
and $\tilde{h}$ form $S_g^0$ to itself representing $t$ and $h$
respectively such that ${\tilde{t}}^k=1$ and
$\tilde{h}\tilde{t}{\tilde{h}}^{-1}\in\<\tilde{t}\>$.}\vspace{1.5mm}

The two theorems above imply that if $\alpha\in\Ga_g^0$ is of finite
order and it point-wise fixes all, say $n$, of its branch points,
then the normalizer of $\alpha$ in $\Ga_g^0$ is isomorphic to that
of $\alpha$ in $\Ga_g^n$ with branch points being the punctures. The
isomorphism between the two normalizers is established by selecting
homeomorphisms of surfaces using Theorem~\ref{BH3} and map them to
each other and Theorem~\ref{BH1} ensures that the maps are well
defined at the mapping class level. More generally, if we remove, or
mark, $m$ branch points before hand, the same argument implies the
following.

{\prop\label{normalizer} Let $m\geqslant 0$, $g\geqslant 0$ be
integers such that $(m,g)\ne(0,0),(0,1)$. Let $\alpha\in P\Ga_g^m$
be of finite order such that there are exactly $n$ branch points in
$S_g^m$ which are point-wise fixed by $\alpha$. By removing these
$n$ additional points from $S_g^m$, one has
$\alpha\in\Ga_g^{n+m,\Sigma_n\times\{1\}^m}$ and the normalizer of
$\alpha$ in $P\Ga_g^m$ is isomorphic to that of $\alpha$ in
$\Ga_g^{n+m,\Sigma_n\times\{1\}^m}$.}\vspace{1.5mm}

Proposition~\ref{normalizer} transfers the problem of finding the
normalizer of a cyclic branched covering map to that of a cyclic
regular covering map. Using this regular covering map, we can obtain
a more detailed description of these normalizers when the quotient
surface is homeomorphic to a punctured 2-sphere. This was partially
done in \cite{BH} and \cite{X1} by showing that there is an
imbedding of the reduced normalizer of $\<\alpha\>$, i.e. the
normalizer of $\<\alpha\>$ modulo $\<\alpha\>$, into the mapping
class group of the quotient surface. We will describe the normalizer
of $\<\alpha\>$ in terms of the automorphism group of the
fundamental group of the quotient surface.

Assume that $\<\alpha\>\cong\Z/p\Z$ is a subgroup of
$\Ga_g^{n+m,\Sigma_n\times\{1\}^m}$ with no additional fixed points
and $S_g^{n+m}/\<\alpha\>$ is homeomorphic to $S_h^{n+m}$, where the
numbers $g$, $h$, $n+m$ and $p$ must satisfy the Riemann-Hurwitz
formula $2g-2=p(2h-2)+(n+m)(p-1)$. Let
$f:\pi_1(S_h^{n+m})\ra\<\alpha\>$ be the group homomorphism obtained
from the regular covering $S^{n+m}_g\ra S_h^{n+m}$ and
$x_1,x_2,\cdots,x_{n+m}$ be the elements in $\pi_1(S_h^{m+n})$ that
represents the $n+m$ punctures. Then the values $f(x_i)$,
$i=1,2,\cdots,n+m$, satisfy $f(x_1)+f(x_1)+\cdots+f(x_{n+m})=0$ and
$f(x_i)\ne 0$ for all $1\leqslant i\leqslant n+m$. We denote by
$(\Z/p\Z)^*$ the non-zero elements in $\Z/p\Z$ and call the vector
$F_\alpha=(f(x_1),f(x_2),\cdots,f(x_{n+m}))\in((\Z/p\Z)^*)^{n+m}$ the
fixed point data vector of $\alpha$. Given a permutation
$\sigma\in\Sigma_{n+m}$ and an element $k\in(\Z/p\Z)^*$ define
\begin{eqnarray*}\sigma(F_\alpha)&=&(f(x_{\sigma(1)}),f(x_{\sigma(2)}),
\cdots,f(x_{\sigma(n+m)}))\in((\Z/p\Z)^*)^{n+m} \mbox{ and }\\
 kF_\alpha&=&(kf(x_1),kf(x_2),\cdots, kf(x_{n+m}))\in((\Z/p\Z)^*)^{n+m}.
\end{eqnarray*} Let
\begin{eqnarray*}\Sigma_\alpha^N&=&\{\gamma\in\Sigma_n\times\{1\}^m|
\gamma(F_\alpha)=kF_\alpha \mbox{ for some }k\in(\Z/p\Z)^*\},\\
\Sigma_\alpha^C&=&\{\gamma\in\Sigma_n\times\{1\}^m|\gamma(F_\alpha)=F_\alpha\}.
\end{eqnarray*} We call $\Sigma_\alpha^N$ the normalizing
permutation group of $\alpha$ and $\Sigma_\alpha^C$ the centralizing
permutation group of $\alpha$. For each $t\in(\Z/p\Z)^*$, let
$l_t=|\{1\leqslant i\leqslant n~|~f(x_i)=t\}|$. We observe from the
definition that
$\Sigma_\alpha^C=\Sigma_{l_1}\times\Sigma_{l_2}
\times\cdots\times\Sigma_{l_{p-1}}\times\{1\}^m$.
Consider the following commutative diagram
$$\minCDarrowwidth0.24in\begin{CD}
1 @>>> \pi_1(S_g^{n+m}) @>i>> \pi_1(S_h^{n+m}) @>f>>\alpha@>>>1\\
 @. @V\sigma VV @V\sigma VV @VVV  \\
1 @>>> {\rm Inn}(\pi_1(S_g^{n+m})) @>i>>{\rm
Aut}_+^{\Sigma_n\times\{1\}^m}(\pi_1(S_g^{n+m}))
@>>>\Ga_g^{\Sigma_n\times\{1\}^m}@>>>1 \\
\end{CD}$$ where the imbedding $\sigma$ (as long as $S_g^{n+m}$
and $S_h^{n+m}$ are not $S_0^1$, $S_0^0$ or $S^0_1$) is given by the
conjugation action on $\pi_1(S_g^{m+n})$, i.e., for every
$u\in\pi_1(S_h^{n+m})$, $\sigma(u)(y)=i^{-1}(ui(y)u^{-1})$ for all
$y\in\pi_1(S_g^{n+m})$. We know that the normalizer of $\<\alpha\>$
in $\Ga_g^{\Sigma_n\times\{1\}^m}$ is isomorphic to that of
$\sigma(\pi_1(S_h^{n+m}))$ in ${\rm
Aut}_+^{\Sigma_n\times\{1\}^m}(\pi_1(S_g^{n+m}))$  modulo $i({\rm
Inn}(\pi_1(S_g^{n+m})))$. {\lem\label{N} The normalizer of
$\sigma(\pi_1(S_h^{n+m}))$ in ${\rm
Aut}_+^{\Sigma_n\times\{1\}^m}(\pi_1(S_g^{n+m}))$ is isomorphic to
${\rm Aut}_+^{\Sigma_\alpha^N}(\pi_1(S_h^{n+m}))$ when $h=0$.}
\vspace{1.5mm}

\proof Let $$K=\{\gamma\in{\rm
Aut}_+^{\Sigma_n\times\{1\}^m}(\pi_1(S_h^{n+m}))~|~
\gamma(i(\pi_1(S_g^{n+m})))=i(\pi_1(S_g^{n+m}))\}$$ and $J$ be the
normalizer of $\sigma(\pi_1(S_h^{n+m}))$ in ${\rm
Aut}_+^{\Sigma_n\times\{1\}^m}(\pi_1(S_g^{n+m}))$. We define two
homomorphisms $\Phi:J\ra K$ and $\Psi:K\ra J$:
\begin{eqnarray*} &&\mbox{for every }x\in J,
\Phi(x)(u)=\sigma^{-1}(x\sigma(u)x^{-1})\mbox{ for all }
u\in\pi_1(S_h^{n+m});\\
&&\mbox{for every }y\in K, \Psi(y)(v)=i^{-1}(y(i(v)))\mbox{ for all
} v\in\pi_1(S_g^{n+m}).\end{eqnarray*} It can be routinely checked
that $\Phi$ and $\Psi$ are inverse to each other. Now we need to
show that $K={\rm Aut}_+^{\Sigma^N_\alpha}(\pi_1(S_h^{n+m}))$ when
$h=0$. Recall that
\begin{eqnarray*} \pi_1(S_0^{n+m})=\langle
x_1,x_2,\cdots,x_{n+m}~|~x_1x_2\cdots x_{n+m}=1\rangle,
\end{eqnarray*} Let $H$ be the normal closure of
$\{t^p~|~t\in\pi_1(S_0^{n+m})\}$ in $\pi_1(S_0^{n+m})$ and
$M=H[\pi_1,\pi_1]$, where $[\pi_1,\pi_1]$ is the commutator subgroup
of $\pi_1(S_0^{n+m})$. Then $M$ is a characteristic subgroup of
$\pi_1(S_0^{n+m})$ and the homomorphism $f$ in the commutative diagram factors through
$$\pi_1(S_0^{n+m})/M=\bigoplus_{i=1}^{n+m}\Z/p\Z x_i
/\<x_1+ x_2+\cdots+x_{n+m}\>.$$ An automorphism in ${\rm
Aut}_+^{\Sigma_n\times\{1\}^m}(\pi_1(S_0^{n+m}))$ fixes the subgroup
$i(\pi_1(S_g^{n+m}))={\rm Ker}(f)$ if and only if its induced
automorphism in $\Aut(\pi_1(S_0^{n+m})/M)$ fixes the kernel of
$$\pi_1(S_0^{n+m})/M\ra\<\alpha\>,$$ which is given by
\begin{eqnarray*}\{t_1x_1+t_2x_2+\cdots+t_{n+m}x_{n+m}\in\bigoplus_{i=1}^{n+m}\Z/p\Z
x_i /\<x_1+
x_2+\cdots+x_{n+m}\>&&\\
|~t_1f(x_1)+t_2f(x_2)+\cdots+t_{n+m}f(x_{n+m})=0\}.&&\end{eqnarray*}
In other words the kernel is the hyperplane $F_\alpha^{\perp}$ in
$\pi_1(S_0^{n+m})/M$ that is orthogonal to the fixed point data
vector $F_{\alpha}$ of $\alpha$. On the other hand, every  automorphism in the
group ${\rm Aut}_+^{\Sigma_n\times\{1\}^m}(\pi_1(S_0^{n+m}))$
induces a permutation automorphism of $\pi_1(S_0^{n+m})/M$ which
permutes the elements $x_1,x_2,\cdots,x_{n}$, and every such
permutation automorphism in ${\Sigma_n\times\{1\}^m}$ is indeed
induced by some automorphism in ${\rm
Aut}_+^{\Sigma_n\times\{1\}^m}(\pi_1(S_0^{n+m}))$. If
$\gamma\in{\Sigma_n\times\{1\}^m}$, one has
$\gamma(F_\alpha^{\perp})=\gamma(F_\alpha)^{\perp}$. Hence
$\gamma(F_\alpha^{\perp})=F_\alpha^{\perp}$ if and only if
$\gamma(F_\alpha)^{\perp}=F_\alpha^{\perp}$, or equivalently
$\gamma\in\Sigma_\alpha^N$ and $K={\rm
Aut}_+^{\Sigma_\alpha^N}(\pi_1(S_0^{n+m}))$, when $h=0$.\qed

{\rem  The group $K$ in the proof is the group that maps loops that
lift to loops to loops that lift to loops as all such loops form the
group $i(\pi_1(S_g^{n+m}))$. This property was used in
\cite{BH,X1}.}\vspace{1.5mm}

Combining Proposition~\ref{normalizer} and Lemma~\ref{N}, we have
proved {\thm\label{NC} Let $m\geqslant 0$, $g\geqslant 0$ be
integers such that $(m,g)\ne(0,0),(0,1)$. Let $\alpha\in P\Ga_g^m$
be an element of prime order $p$ with exactly $n$ branch points in
$S_g^m$. Let $\Sigma_\alpha^N$ and $\Sigma_\alpha^C$ be the
normalizing and centralizing permutation groups of $\alpha$
respectively. If the quotient surface $S_g^{m}/\<\alpha\>$ is
homeomorphic to $S_0^{m}$, then the normalizer of $\alpha$ in
$P\Ga_g^m$ is isomorphic to ${\rm
Aut}_+^{\Sigma_\alpha^N}(\pi_1(S_0^{n+m}))/{\rm
Inn}(\pi_1(S_g^{n+m}))$ and the centralizer of $\alpha$ in
$P\Ga_g^m$ is isomorphic to ${\rm
Aut}_+^{\Sigma_\alpha^C}(\pi_1(S_0^{n+m}))/{\rm
Inn}(\pi_1(S_g^{n+m}))$.}

\section{The $p$-primary part of the Farrell cohomology of certain pure mapping
class groups} The mapping class groups are known to have finite
virtual cohomological dimensions (vcd)(see \cite{Lu2,M}), i.e. they
have subgroups of finite index which have finite cohomological
dimensions. For group $\Ga$ of finite vcd, the ordinary cohomology
and Farrell cohomology are isomorphic above the vcd, i.e. $$H^k(\Ga,
M)\cong\widehat{H}^i(\Ga,M)$$ for all $\Ga$-module $M$ and all
$i>{\rm vcd}(\Ga)$.  Such a group $\Ga$ is said to have $p$-periodic
cohomology for some prime $p$ if there is a positive integer $d$
such that
$\widehat{H}^i(\Ga,M)_{(p)}\cong\widehat{H}^{i+d}(\Ga,M)_{(p)}$ for
all integer $i$ and all $\Ga$-module $M$, where
$\widehat{H}^i(\Ga,M)_{(p)}$ is the $p$-primary part of
$\widehat{H}^i(\Ga,M)$. By \cite[Theorem 6.7]{Br}, a group of finite
vcd is $p$-periodic for a prime $p$ if and only if $\Gamma$ is of
$p$-rank 1, i.e. every finite elementary abelian $p$-subgroup is
rank 1. The following theorem is needed in our study of the
cohomology groups of certain pure mapping class groups.
{\thm\label{Brown}\cite[Corollary 7.4]{Br} Suppose that $\Gamma$ is
of finite vcd and every elementary abelian $p$-subgroup of $\Gamma$
has rank $\leqslant 1.$ Then
\begin{eqnarray*}
{\widehat H}^i(\Gamma,{\mathbb
Z})_{(p)}\cong\prod_{P\in\wp}{\widehat H}^i(N_{\Gamma}(P) ,{\mathbb
Z})_{(p)},
\end{eqnarray*}
where $\wp$ is a set of representatives for the conjugacy classes of
subgroups of $\Gamma$ of order $p$, $N_{\Gamma}(P)$ is the
normalizer of $P$ in $\Gamma$.}\vspace{1.5mm}

In this section we study the $p$-primary part of the Farrell
cohomology of the pure mapping class groups. Lu \cite{Lu1,Lu2,Lu3}
determined almost all of the $p$-primary part of Farrell cohomology
of pure mapping class groups of low genus, as well as
$P\Gamma_{n(p-1)/2}^{m}$ for small $n$. We obtain some general
results on the $p$-primary part of the Farrell cohomology of
$P\Gamma_{n(p-1)/2}^{m}$ for some $m$'s, and the 2-primary part of the Farrell
cohomology of the pure mapping class group $P\Gamma_{n}^{m}$ for
$m=2n-1,2n,2n+1,2n+2$, where $p$ is an odd prime and $n\geqslant 1$
is an integer. The explicit calculation is carried out for
${\widehat H}^i(P\Gamma_{n(p-1)/2}^{n+1},{\mathbb Z})_{(p)}$,
${\widehat H}^i(P\Gamma_{n(p-1)/2}^{n+2},{\mathbb Z})_{(p)}$,
${\widehat H}^i(P\Gamma_{n}^{2n+1},{\mathbb Z})_{(2)}$ and
${\widehat H}^i(P\Gamma_{n}^{2n+2},{\mathbb Z})_{(2)}$.  The reason
why we choose these groups is because the Riemann-Hurwitz formula
implies {\lem\label{reason}
 Let $p$ be any odd prime, $m\geqslant 0$,  and $n\geqslant 1$ be two integers.
\begin{itemize}
\item[(1)] If $n-p+3\leqslant m\leqslant n+2$,
then every subgroup
$\<\alpha\>$ of order $p$ in $P\Ga^m_{n(p-1)/2}$ has $n+2$ fixed
points including the punctures and the quotient surface
$S^m_{n(p-1)/2}/\<\alpha\>$ is homeomorphic to $S^m_0$.
\item[(2)] If $2n-1\leqslant m\leqslant 2n+2$, then every subgroup
$\<\alpha\>$ of order $2$ in $P\Ga^{m}_{n}$ has $2n+2$ fixed points
including the punctures and the quotient surface
$S^m_{n}/\<\alpha\>$ is homeomorphic to
$S^m_0$.\end{itemize}}

{\rem Lemma 5.2  is also a consequence of the main result of
\cite{GJ}: Every non-free orientable action of $\Z/p\Z$ on
$S^{0}_{g}$, $p$ odd, is given by the connected sum at fixed points
of the canonical actions of $\Z/p\Z$ on $S^0_p$ and $S^0_{(p-1)/2}$.
For $p=2$ the canonical action is for $\Z/2\Z$ on $S^0_1$ and
$S^0_2$.}

In view of Theorem~\ref{Brown}, we start our computation by
obtaining some structural information about the cohomology groups of
normalizers of elements of prime order in $P\Ga_g^m$ such that the
quotient surface is homeomorphic to $S_0^m$.

{\lem\label{ele} Let $\alpha\in P\Ga^m_g$ be an element of prime
order $p$ such that $S_g^m/\<\alpha\>$ is homeomorphic to $S^m_0$,
and $\Sigma_\alpha^C$ its centralizing permutation group. Let
$N(\alpha)$ be the normalizer of $\alpha$ in $P\Ga^m_g$. If
$\Sigma_\alpha^C$ contains a subgroup of the form $K=H\times\{1\}^l$
for some $l\geqslant 3$, i.e. $K$ fixes at least three punctures,
and the index $[\Sigma_\alpha^C:K]$ is relatively prime to $p$, then
$\widehat{H}^*(N(\alpha),\Z)_{(p)}$ is an elementary abelian
group.}\vspace{1.5mm}

\proof By Theorem~\ref{NC}, the centralizer $C(\alpha)$ of $\alpha$
in $P\Ga^m_g$ is isomorphic to the group ${\rm
Aut}_+^{\Sigma_\alpha^C}(\pi_1(S_0^{n+m}))/{\rm
Inn}(\pi_1(S_g^{n+m}))$, which will also be called $C(\alpha)$.  Let
$K(\alpha)={\rm Aut}_+^K(\pi_1(S_0^{n+m}))/{\rm
Inn}(\pi_1(S_g^{n+m}))$. Then
$[C(\alpha):K(\alpha)]=[\Sigma_\alpha^C:K]$ is prime to $p$. It
follows that
$[N(\alpha):K(\alpha)]=[N(\alpha):C(\alpha)][C(\alpha):K(\alpha)]$
is also prime to $p$. Since the group $K=H\times\{1\}^l$ and
$l\geqslant 3$, by the algebraic description of the Artin braid group
and the mapping class group and Theorem~\ref{split},
$${\rm Aut}_+^K(\pi_1(S_0^{n+m}))=\<B_{n+m-1}^{H\times\{1\}^{l-1}},
{\rm Inn}(\pi_1(S_0^{m+n}))\>=\<\Ga_0^{n+m,K}\times\Z, {\rm
Inn}(\pi_1(S_0^{m+n})) \>,$$ where the center
$\Z=B_{n+m-1}^{H\times\{1\}^{l-1}}\cap{\rm Inn}(\pi_1(S_0^{m+n}))$,
and $\<A,B\>$ denotes the subgroup generated by $A$ and $B$.
Therefore,
\begin{align*}K(\alpha)=~&{\rm
Aut}_+^K(\pi_1(S_0^{n+m}))/{\rm
Inn}(\pi_1(S_g^{m+n}))\\=~&\<B_{n+m-1}^{H\times\{1\}^{l-1}}, {\rm
Inn}(\pi_1(S_0^{m+n}))\>/{\rm Inn}(\pi_1(S_g^{m+n}))
\end{align*} is generated by the group $B_{n+m-1}^{H\times\{1\}^{l-1}}/
[B_{n+m-1}^{H\times\{1\}^{l-1}}\cap{\rm Inn}(\pi_1(S_g^{m+n}))]$ and
the group ${\rm Inn}(\pi_1(S_0^{m+n}))/{\rm
Inn}(\pi_1(S_g^{m+n}))=\<\alpha\>$. Since
$$B_{n+m-1}^{H\times\{1\}^{l-1}}=\Ga_0^{n+m,K}\times(B_{n+m-1}^{H\times\{1\}^{l-1}}\cap{\rm
Inn}(\pi_1(S_0^{m+n}))),$$ one has
\begin{align*}B_{n+m-1}^{H\times\{1\}^{l-1}}&/
[B_{n+m-1}^{H\times\{1\}^{l-1}}\cap{\rm Inn}(\pi_1(S_g^{m+n}))]=\\
&\Ga_0^{n+m,K}\times[(B_{n+m-1}^{H\times\{1\}^{l-1}}\cap{\rm
Inn}(\pi_1(S_0^{m+n})))/(B_{n+m-1}^{H\times\{1\}^{l-1}}\cap{\rm
Inn}(\pi_1(S_g^{m+n})))].
\end{align*}
The group $(B_{n+m-1}^{H\times\{1\}^{l-1}}\cap{\rm
Inn}(\pi_1(S_0^{m+n})))/(B_{n+m-1}^{H\times\{1\}^{l-1}}\cap{\rm
Inn}(\pi_1(S_g^{m+n})))$ is a subgroup of $\<\alpha\>={\rm
Inn}(\pi_1(S_0^{m+n}))/{\rm Inn}(\pi_1(S_g^{m+n}))\>$ and therefore
\begin{align*}K(\alpha)=\Ga_0^{n+m,K}\times\<\alpha\>.
\end{align*}  By the K\"unneth formula, $H^*(K(\alpha),\Z)$ is an
elementary abelian group in sufficiently high dimension as
$\Ga_0^{n+m,K}$ has finite cohomological dimension. Since
$K(\alpha)$ is $p$-periodic as it has $p$-rank 1,
$\widehat{H}^*(K(\alpha),\Z)$ is elementary abelian. According to
\cite[Proposition 9.5]{Br}, the composition of restriction map ${\rm
res}^{N(\alpha)}_{K(\alpha)}:\widehat{H}^*(N(\alpha),\Z)\ra
\widehat{H}^*(K(\alpha),\Z)$ with the transfer map ${\rm
cor}_{N(\alpha)}^{K(\alpha)}:\widehat{H}^*(K(\alpha),\Z)\ra
\widehat{H}^*(N(\alpha),\Z)$ yields the scalar multiplication map
multiplying with $[N(\alpha):K(\alpha)]$, i.e. $${\rm
cor}_{N(\alpha)}^{K(\alpha)}{\rm
res}^{N(\alpha)}_{K(\alpha)}(z)=[N(\alpha):K(\alpha)]z$$ for all
$z\in\widehat{H}^*(N(\alpha),\Z)$. Hence
$\widehat{H}^*(N(\alpha),\Z)$ is annihilated by $p[N(\alpha):K(\alpha)]$
as ${\rm res}^{N(\alpha)}_{K(\alpha)}(pz)=p({\rm
res}^{N(\alpha)}_{K(\alpha)}(z))=0$, and
$\widehat{H}^*(N(\alpha),\Z)_{(p)}$ is an elementary abelian
group.\qed

{\cor\label{co} Let $p$ be an odd prime and $\alpha\in
P\Ga_{n(p-1)/2}^m$ an element of order $p$ with $n+2$ fixed points
including punctures. Let $N(\alpha)$ be the normalizer of $\alpha$
in $P\Ga_{n(p-1)/2}^m$. If $m$ and $n$ satisfy
\begin{itemize}
\item[(i)] $m\geqslant 3$; or
\item[(ii)] $1\leqslant m\leqslant 2$ and $n\not\equiv 0~(\!\!\!\!\mod p)$; or
\item[(iii)] $m=0$ and $n\not\equiv 0$ or $-2~(\!\!\!\!\mod p)$,
\end{itemize} then $\widehat{H}^*(N(\alpha),\Z)_{(p)}$ is an elementary
abelian group.}\vspace{1.5mm}

\proof By the Riemann-Hurwitz formula, the quotient surface
$S_{n(p-1)/2}^m/\<\alpha\>$ is homeomorphic to $S_{0}^m$. The case
$m\geqslant 3$ follows rather obviously from Lemma ~\ref{ele}. If
$m=0$ and $n\not\equiv 0$ or $-2 (\!\!\!\mod p)$, let
$\Sigma_{\alpha}^C$ be the centralizing permutation group of
$\alpha$. Then
$\Sigma_{\alpha}^C=\Sigma_{l_1}\times\Sigma_{l_2}\times\cdots\times\Sigma_{l_{p-1}}$
with \begin{align*}&l_1+l_2+\cdots+l_{p-1}=n+2\not\equiv
0~(\!\!\!\!\mod p)\\ &l_1+2l_2+\cdots+(p-1)l_{p-1}\equiv 0~
(\!\!\!\!\mod p).\end{align*} If there are three $l_i$'s, say $l_1$,
$l_2$ and $l_3$, prime to $p$, we can choose
$H=\Sigma_{l_1-1}\times\Sigma_{l_2-1}\times\Sigma_{l_3-1}\times
\Sigma_{l_4}\cdots\times\Sigma_{l_{p-1}}$. If there are only two
$l_i$'s, say $l_1$ and $l_2$, prime to $p$, then
$(l_1-1)+(l_2-1)\equiv n\not\equiv 0$ $(\!\!\!\!\mod p)$ and there
is at least one, say $l_1-1$, prime to $p$. Hence we can choose
$H=\Sigma_{l_1-2}\times\Sigma_{l_2-1}\times\Sigma_{l_3}\times
\cdots\times\Sigma_{l_{p-1}}$. By Lemma~\ref{ele},
$\widehat{H}^*(N(\alpha),\Z)_{(p)}$ is elementary abelian. The case
$1\leqslant m\leqslant 2$ and $n\not\equiv 0~(\!\!\!\!\mod p)$ can
be similarly proved.\qed \vspace{1mm}

As we can see in most of the situations with quotient being the
punctured 2-sphere, the $p$-primary part of Farrell cohomology has
exponent $p$. By Lemma~\ref{reason} and Corolary~\ref{co}, we have

{\thm\label{ppart}  Let $m\geqslant 0$ and $n\geqslant 1$ be two
integers.
\begin{itemize}
\item[(1)] If  $n-p+3\leqslant m\leqslant n+2$
then $\widehat{H}^*(P\Ga^m_{n(p-1)/2}, \Z)_{(p)}$ is an elementary
abelian group.
\item[(2)] If $2n-1\leqslant m\leqslant 2n+2$, then $\widehat{H}^*(P\Ga^m_{n},
\Z)_{(2)}$ is an elementary abelian group when $P\Ga_n^m\ne
P\Ga_1^1$ or $P\Ga_1^2$.\end{itemize} } \vspace{1.5mm}

We are now ready to compute the cohomology groups
$\widehat{H}^i(P\Gamma_{n(p-1)/2}^{n+1},{\mathbb Z})_{(p)}$,
$\widehat {H}^i(P\Gamma_{n(p-1)/2}^{n+2},{\mathbb Z})_{(p)}$,
$\widehat {H}^i(P\Gamma_{n}^{2n+1},{\mathbb Z})_{(2)}$ and $\widehat
{H}^i(P\Gamma_{n}^{2n+2},{\mathbb Z})_{(2)}$. By
Theorem~\ref{Brown}, we need to find all conjugacy classes of
subgroups of order $p$, which was done in \cite{Lu1,Lu2}. The idea
is that there is an one to one correspondence between the conjugacy
classes of subgroups $\<\alpha\>$ of order $p$ in
$P\Gamma_{n(p-1)/2}^{n+1}$, as well as in
$P\Gamma_{n(p-1)/2}^{n+2}$, and the fixed point data vector up to
scalers $\{kF_\alpha~|~k\in\Z/p\Z\}$. Hence the number of conjugacy
classes of subgroups of order $p$ in $P\Gamma_{n(p-1)/2}^{n+1}$ and
in $P\Gamma_{n(p-1)/2}^{n+2}$ is equal to the number of solutions to
the equation $f(x_1)+f(x_2)+\cdots+f(x_{n+1})=1$ in $\Z/p\Z$ and
$f(x_i)\ne 0$ for all $i=1,2,\cdots,n+1$. If we denote the number of
such solutions by $\nu(n+1)$, then we can get the following
recursion formula based upon whether or not $f(x_{n})=1$,
$$(p-1)\nu(n-1)+(p-2)\nu(n)=\nu(n+1).$$  The initial values are
obviously $\nu(1)=1$ and $\nu(2)=p-2$, By solving this recursion,
one gets $\nu(n+1)=[(p-1)^{n+1}+(-1)^n]/p$.  The groups
$P\Ga_n^{2n+1}$ and $P\Ga_n^{2n+2}$ obviously have only one
conjugacy class of subgroups of order 2 which fixes $2n+2$ points.
As a result,
{\lem\label{num} The groups  $P\Ga_{n(p-1)/2}^{n+1}$ and
$P\Ga_{n(p-1)/2}^{n+2}$ have $[(p-1)^{n+1}+(-1)^n]/p$ conjugacy
classes of subgroups of order $p$, and the groups $ P\Ga_{n}^{2n+1}$
and $P\Ga_{n}^{2n+2}$ have $1$ conjugacy class of subgroups of order
$2$.}\vspace{1.5mm}

If one follows the same line of calculation as in Lemma~\ref{ele},
then for any $\alpha\in P\Ga_{n(p-1)/2}^{n+1}$ or $\alpha\in
P\Ga_{n(p-1)/2}^{n+2}$ of order $p$, one can see that
$$N(\alpha)=C(\alpha)=K(\alpha)=P\Ga_0^{n+2}\times\Z/p\Z.$$
For any $\alpha\in P\Ga_{n}^{2n+1}$ or $\alpha\in P\Ga_{n}^{2n+2}$
of order 2, one has
$$N(\alpha)=C(\alpha)=K(\alpha)=P\Ga_0^{2n+2}\times\Z/2\Z.$$

The cohomology groups of $P\Gamma_0^m$ have been computed by Cohen
\cite{Co}. These are free abelian groups. The Poincar\'e series of
the cohomology is given by
$$P_m(t)=(1+2t)(1+3t)\cdots(1+(m-2)t).$$ Using the K\"unneth
formula, one finds that
\begin{align*}
 H^i(N(\alpha),{\mathbb Z})_{(p)}
&=({\mathbb Z}/p{\mathbb Z})^{(P_{n+2}(1)+(-1)^iP_{n+2}(-1))/2}\\
&=({\mathbb Z}/p{\mathbb
Z})^{(\frac{1}{2}(n+1)!+(-1)^{i+n-1}(n-1)!)/2}\\
&=({\mathbb Z}/p{\mathbb Z})^{(n-1)![n^2+n-(-1)^{n+i}2]/4}
\end{align*}
for $i$ sufficiently large and $\alpha\in P\Ga_{n(p-1)/2}^{n+1}$ or
$P\Ga_{n(p-1)/2}^{n+2}$, and
\begin{align*}
 H^i(N(\alpha),{\mathbb Z})_{(2)}
&=({\mathbb Z}/2{\mathbb Z})^{(P_{2n+2}(1)+(-1)^iP_{2n+2}(-1))/2}\\
&=({\mathbb Z}/2{\mathbb
Z})^{(\frac{1}{2}(2n+1)!+(-1)^{i+2n-1}(2n-1)!)/2}\\
&=({\mathbb Z}/2{\mathbb Z})^{(2n-1)![2n^2+n-(-1)^{i}]/2}
\end{align*}
for $i$ sufficiently large and $\alpha\in P\Ga_{n}^{2n+1}$ or
$P\Ga_{n}^{n+2}$. Since pure mapping class groups are of finite vcd,
we know that their Farrell cohomology groups coincide with their
ordinary cohomology groups above the vcd. Since pure mapping class
groups are periodic, by Theorem~\ref{Brown} and Lemma~\ref{num}, one
has {\thm\label{explicit} For any given odd prime $p$ and positive
integer $n$,   the $p$-primary part of the Farrell cohomology of
$P\Gamma_{n(p-1)/2}^{n+1}$ and $P\Gamma_{n(p-1)/2}^{n+2}$ is given
by
\begin{eqnarray*}&&{\widehat H}^i(P\Gamma_{n(p-1)/2}^{n+1},{\mathbb
Z})_{(p)} \cong{\widehat H}^i(P\Gamma_{n(p-1)/2}^{n+2},{\mathbb
Z})_{(p)}\\&&\cong({\mathbb Z}/ p{\mathbb
Z})^{(n-1)![(p-1)^{n+1}+(-1)^n][n^2+n-(-1)^{n+i}2]/(4p)}.\end{eqnarray*}
The $2$-primary part of the Farrell cohomology of $P\Gamma_{n}^{2n+1}$
and $P\Gamma_{n}^{2n+2}$ is given by
\begin{eqnarray*}&&{\widehat H}^i(P\Gamma_{n}^{2n+1},{\mathbb
Z})_{(2)} \cong{\widehat H}^i(P\Gamma_{n}^{2n+2},{\mathbb
Z})_{(2)}\cong({\mathbb Z}/ 2{\mathbb
Z})^{(2n-1)![2n^2+n-(-1)^{i}]/2}.\end{eqnarray*}}

As illustrations of Theorem 5.8, we conclude this paper with the following remark.

{\rem cf.(\cite{Lu2,Lu3})
\begin{align*}
&{\widehat H}^0(P\Gamma^3_1,\Z)_{(2)}\cong\Z/2\Z,
~{\widehat H}^1(P\Gamma^3_1,\Z)_{(2)}\cong(\Z/2\Z)^2,\\
&{\widehat H}^0(P\Gamma^2_{(p-1)/2},\Z)_{(p)}\cong{\widehat H}^0(P\Gamma^3_{(p-1)/2},\Z)_{(p)}\cong(\Z/p\Z)^{p-2},\\
&{\widehat H}^1(P\Gamma^2_{(p-1)/2},\Z)_{(p)}\cong{\widehat H}^1(P\Gamma^3_{(p-1)/2},\Z)_{(p)}\cong  0,\\
&{\widehat H}^0(P\Gamma^3_{p-1},\Z)_{(p)}\cong{\widehat H}^0(P\Gamma^4_{p-1},\Z)_{(p)}\cong(\Z/p\Z)^{p^2-3p+3},\\
&{\widehat H}^1(P\Gamma^3_{p-1},\Z)_{(p)}\cong{\widehat H}^1(P\Gamma^4_{p-1},\Z)_{(p)}\cong (\Z/p\Z)^{2(p^2-3p+3)},\\
&{\widehat H}^0(P\Gamma^4_{3(p-1)/2},\Z)_{(p)}\cong{\widehat H}^0(P\Gamma^5_{3(p-1)/2},\Z)_{(p)}\cong (\Z/p\Z)^{7(p^3-4p^2+6p-4)},\\
&{\widehat H}^1(P\Gamma^4_{3(p-1)/2},\Z)_{(p)}\cong{\widehat H}^1(P\Gamma^5_{3(p-1)/2},\Z)_{(p)}\cong (\Z/p\Z)^{5(p^3-4p^2+6p-4)}.
\end{align*} }

\end{document}